\documentclass[12pt]{article}
\usepackage{amsmath}
\usepackage{amssymb,amsfonts,amsmath,latexsym,cite,color, psfrag,graphicx,ifpdf,ytableau,pgfkeys,pgfopts,xcolor,extarrows}

\newtheorem{theo}{Theorem}[section]
\newtheorem{defi}[theo]{Definition}

\newtheorem{prop}[theo]{Proposition}

\allowdisplaybreaks

\makeatletter \@addtoreset{equation}{section}
\@addtoreset{theo}{}\makeatother

\usepackage{hyperref}

\def\qed{\hfill \rule{4pt}{7pt}}
\def\pf{\noindent {\it Proof.} }

\begin{document}
\begin{center}

{\Large\bf Standard Rothe Tableaux }\\[12pt]
\begin{center}
{ Neil J. Y. Fan  }

\vskip 4mm
Department of Mathematics\\
Sichuan University, Chengdu, Sichuan 610064, P.R. China

\vskip 4mm
 fan@scu.edu.cn

\end{center}

\end{center}

\vskip 6mm

\begin{abstract}
Edelman and Greene constructed a bijection between the set of standard Young tableaux and the set of  balanced Young tableaux of the same shape. Fomin,  Greene,   Reiner and  Shimozono introduced the notion of balanced Rothe tableaux of a permutation $w$, and established a bijection between the set of balanced Rothe tableaux of $w$ and the set of reduced words of $w$.
In this paper, we introduce the notion of standard Rothe tableaux of $w$, which are tableaux obtained by labelling  the cells of the Rothe diagram of $w$ such that each row and each column is  increasing.  We show that the number of standard Rothe tableaux of $w$ is smaller than or equal to the number of balanced Rothe tableaux of $w$, with equality if and only if $w$ avoids the four patterns 2413, 2431, 3142 and 4132.  When $w$ is a dominant permutation, i.e., 132-avoiding, the Rothe diagram of $w$ is a Young diagram, so this reduces to the result  of  Edelman and Greene.
\end{abstract}


\section{Introduction}

Let $S_n$ denote the group of permutations on $\{1,2,\ldots,n\}$.
Given $w_0=n\cdots21\in S_n$, Stanley \cite{Stanley} showed algebraically that the number of reduced words of $w_0$ is equal to the number of standard staircase Young tableaux of shape $(n-1,\ldots,1)$.
In  \cite{Edelman}, Edelman and Greene introduced the notion of balanced Young tableaux, and showed that the number of balanced staircase Young tableaux of shape $(n-1,\ldots,1)$ is equal to the number of reduced words of $w_0$. Furthermore, they constructed a bijection between the set of standard Young tableaux and the set of balanced Young tableaux of any shape, which proved the result of Stanley combinatorially.

Notice that the Rothe diagram  of a dominant permutation, i.e., 132-avoiding, is a Young diagram, see, for example, \cite{Ma}. It is natural to generalize the results of \cite{Edelman} to the Rothe tableaux of a permutation.  In \cite{Fomin}, Fomin,  Greene,   Reiner and   Shimozono introduced the notion of balanced Rothe tableaux of any permutation  $w$, and constructed a simple bijection   between the set of reduced words of  $w$ and the set of balanced Rothe tableaux of $w$.

In this paper, we introduce the notion of standard Rothe tableaux of  $w$, which are tableaux obtained by labelling the cells of the Rothe diagram   of $w$ with consecutive integers $1,2,\ldots,$ such that each row and each column is increasing.
The   result  of Edelman and Greene \cite{Edelman} can be viewed as a bijection between   standard Rothe  tableaux   and   balanced Rothe tableaux for dominant permutations.
For a general permutation $w$, we show that the number of  standard Rothe tableaux of $w$ is smaller than or equal to the number of  balanced Rothe tableaux of $w$, with equality if and only if $w$ avoids the four patterns 2413, 2431, 3142 and 4132. Moreover, when the equality holds, we provide an explicit formula for such number of standard Rothe tableaux.

Let $SRT(w)$ denote the set of standard Rothe tableaux of $w$. If $w$ avoids the four patterns 2413, 2431, 3142 and 4132, we show that each connected component of the Rothe diagram of $w$ is a Young diagram, and any two Young diagrams are nonintersecting, i.e., have no cells in the same row or column. Thus $|SRT(w)|$ can be expressed as the product of  the numbers of some standard Young tableaux.
 On the other hand, if $w$ contains one of these four patterns, the main steps of our proof are to first establish an injection  $\eta$, called the lifting operation, from $SRT(w)$ to $SRT(\widetilde{w})$, where $\widetilde{w}$ is a dominant permutation.
Then by using the lifting operation $\eta$,  we show that there is an injective but not surjective map from
$SRT(w)$  to the set of reduced words of $w$.

We proceed to recall some notation and definitions.
It is well known that $S_n$ is a Coxeter group with generator set $S=\{s_1,s_2,\ldots,s_{n-1}\}$, where $s_i=(i,i+1)$ is the adjacent transposition interchanging $i$ and $i+1$. Let $w\in S_n$ be a permutation, $w$ can be expressed as a product of generators $w=s_{i_1}s_{i_2}\cdots s_{i_k}$. When the number of generators needed is the smallest, such a word is called reduced, and this number is called the length of $w$, denoted by $\ell(w)$. The set of reduced words of $w$ is denoted by $R(w)$.
Given a permutation $w=w_1w_2\cdots w_n$, we say that $w$ avoids the pattern 132, if there does not exist $i_1<i_2<i_3$ such that $w_{i_1}<w_{i_3}<w_{i_2}$. In this case, we also say that $w$ is a dominant permutation, or $w$ is 132-avoiding. We can define other pattern avoiding permutations similarly.

Recall that the Rothe diagram  $D(w)$ of a  permutation
 $w=w_1w_2\cdots w_n$ is a diagram  characterization of the inversions of $w$.
Consider  an $n\times n$ square grid, where we use $(i,j)$ to denote the cell in row $i$
and column $j$. We put a dot in the cell $(i,w_i)$ for $1\leq i\leq n$.
Then the Rothe diagram $D(w)$  is the collection  of cells $B$ such that there is a dot to the right of $B$   and there is a dot below $B$. A standard Rothe tableau of $w$ is obtained by labelling the cells of $D(w)$ with integers $1,2,\ldots,\ell(w)$ such that each row and each column is strictly increasing.

To each cell $B(i,j)$ of $D(w)$, the hook $H_{i,j}(D)$ is the set of cells $(i',j')$ of $D(w)$ such that $i=i'$ and $j'\ge j$, or $j=j'$ and $i'\ge i$. A labelling of the cells of $D(w)$ is called balanced, if for each hook $H_{i,j}(D)$, after rearranging the labels in $H_{i,j}(D)$ increasingly from right to left and from top to bottom, the label of the cell $(i,j)$ remains unchanged.
A balanced Rothe tableau of $w$ is a tableau obtained by giving a balanced labelling of $D(w)$ with integers $1,2,\ldots,\ell(w)$. Denote the set of balanced Rothe tableaux of $w$ by $BRT(w)$.
In particular, when $D(w)$ is a Young diagram of shape $\lambda$, $BRT(w)$ is also the set of balanced Young tableaux of shape $\lambda$, denoted by $BYT(\lambda)$.

\begin{figure}[h,t]
\setlength{\unitlength}{0.5mm}
\begin{center}
\begin{picture}(210,65)
\qbezier[60](0,0)(30,0)(60,0)\qbezier[60](0,10)(30,10)(60,10)
\qbezier[60](0,20)(30,20)(60,20)\qbezier[60](0,30)(30,30)(60,30)
\qbezier[60](0,40)(30,40)(60,40)\qbezier[60](0,50)(30,50)(60,50)
\qbezier[60](0,60)(30,60)(60,60)

\qbezier[60](0,0)(0,30)(0,60)\qbezier[60](10,0)(10,30)(10,60)
\qbezier[60](20,0)(20,30)(20,60)\qbezier[60](30,0)(30,30)(30,60)
\qbezier[60](40,0)(40,30)(40,60)\qbezier[60](50,0)(50,30)(50,60)
\qbezier[60](60,0)(60,30)(60,60)

\put(55,35){\circle*{3}}\put(45,5){\circle*{3}}
\put(25,25){\circle*{3}}\put(5,15){\circle*{3}}
\put(35,55){\circle*{3}}
\put(15,45){\circle*{3}}

\put(0,30){\line(1,0){10}}\put(0,40){\line(1,0){10}}
\put(0,20){\line(1,0){10}}
\put(0,50){\line(1,0){30}}\put(0,60){\line(1,0){30}}
\put(0,20){\line(0,1){40}}\put(10,20){\line(0,1){40}}
\put(20,50){\line(0,1){10}}\put(30,50){\line(0,1){10}}

\put(20,30){\line(1,0){10}}\put(20,40){\line(1,0){10}}
\put(20,30){\line(0,1){10}}\put(30,30){\line(0,1){10}}

\put(40,30){\line(1,0){10}}\put(40,40){\line(1,0){10}}
\put(40,30){\line(0,1){10}}\put(50,30){\line(0,1){10}}


\qbezier[60](80,0)(110,0)(140,0)\qbezier[60](80,10)(110,10)(140,10)
\qbezier[60](80,20)(110,20)(140,20)\qbezier[60](80,30)(110,30)(140,30)
\qbezier[60](80,40)(110,40)(140,40)\qbezier[60](80,50)(110,50)(140,50)
\qbezier[60](80,60)(110,60)(140,60)

\qbezier[60](80,0)(80,30)(80,60)\qbezier[60](90,0)(90,30)(90,60)
\qbezier[60](100,0)(100,30)(100,60)\qbezier[60](110,0)(110,30)(110,60)
\qbezier[60](120,0)(120,30)(120,60)\qbezier[60](130,0)(130,30)(130,60)
\qbezier[60](140,0)(140,30)(140,60)

\put(135,35){\circle*{3}}\put(125,5){\circle*{3}}
\put(105,25){\circle*{3}}\put(85,15){\circle*{3}}
\put(115,55){\circle*{3}}
\put(95,45){\circle*{3}}

\put(80,30){\line(1,0){10}}\put(80,40){\line(1,0){10}}
\put(80,20){\line(1,0){10}}
\put(80,50){\line(1,0){30}}\put(80,60){\line(1,0){30}}
\put(80,20){\line(0,1){40}}\put(90,20){\line(0,1){40}}
\put(100,50){\line(0,1){10}}\put(110,50){\line(0,1){10}}

\put(100,30){\line(1,0){10}}\put(100,40){\line(1,0){10}}
\put(100,30){\line(0,1){10}}\put(110,30){\line(0,1){10}}

\put(120,30){\line(1,0){10}}\put(120,40){\line(1,0){10}}
\put(120,30){\line(0,1){10}}\put(130,30){\line(0,1){10}}

\put(83,52){1}\put(83,42){2}\put(83,32){4}\put(83,22){5}
\put(93,52){3}\put(103,52){6}

\put(103,32){7}  \put(123,32){8}

\qbezier[60](160,0)(190,0)(220,0)\qbezier[60](160,10)(190,10)(220,10)
\qbezier[60](160,20)(190,20)(220,20)\qbezier[60](160,30)(190,30)(220,30)
\qbezier[60](160,40)(190,40)(220,40)\qbezier[60](160,50)(190,50)(220,50)
\qbezier[60](160,60)(190,60)(220,60)

\qbezier[60](160,0)(160,30)(160,60)\qbezier[60](170,0)(170,30)(170,60)
\qbezier[60](180,0)(180,30)(180,60)\qbezier[60](190,0)(190,30)(190,60)
\qbezier[60](200,0)(200,30)(200,60)\qbezier[60](210,0)(210,30)(210,60)
\qbezier[60](220,0)(220,30)(220,60)

\put(215,35){\circle*{3}}\put(205,5){\circle*{3}}
\put(185,25){\circle*{3}}\put(165,15){\circle*{3}}
\put(195,55){\circle*{3}}
\put(175,45){\circle*{3}}

\put(160,30){\line(1,0){10}}\put(160,40){\line(1,0){10}}
\put(160,20){\line(1,0){10}}
\put(160,50){\line(1,0){30}}\put(160,60){\line(1,0){30}}
\put(160,20){\line(0,1){40}}\put(170,20){\line(0,1){40}}
\put(180,50){\line(0,1){10}}\put(190,50){\line(0,1){10}}

\put(180,30){\line(1,0){10}}\put(180,40){\line(1,0){10}}
\put(180,30){\line(0,1){10}}\put(190,30){\line(0,1){10}}

\put(200,30){\line(1,0){10}}\put(200,40){\line(1,0){10}}
\put(200,30){\line(0,1){10}}\put(210,30){\line(0,1){10}}

\put(163,52){3}\put(163,42){1}\put(163,32){7}\put(163,22){5}
\put(173,52){4}\put(183,52){2}

\put(183,32){8}  \put(203,32){6}
\put(25,-10){{\small (a)}}\put(105,-10){{\small (b)}}\put(185,-10){{\small (c)}}
\end{picture}
\end{center}
\caption{ $D(w)$ for $w=426315$, a standard Rothe tableau and a balanced Rothe tableau of $w$.}\label{fig1}
\end{figure}
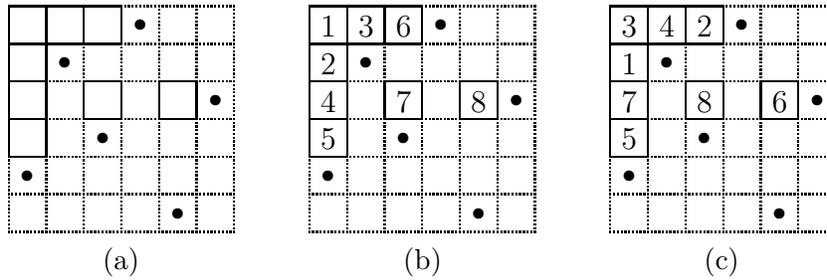

For example, Figure \ref{fig1}(a) gives the Rothe diagram $D(w)$ for $w=426315$,  Figure \ref{fig1}(b) provides a standard Rothe tableau of $w$, and Figure \ref{fig1}(c) illustrates a balanced Rothe tableau of $w$.

\section{Promotion on tableaux}

In this section, we review some basic  definitions and operations on Young tableaux  and recall some results of Edelman and Greene in \cite{Edelman} that we need.

A partition of $n$ is a nonnegative integer sequence $\lambda=(\lambda_1,\ldots,\lambda_k)$, denoted $\lambda\vdash n$, such that $\lambda_1\ge\cdots\ge\lambda_k\ge0$ and $|\lambda|=\lambda_1+\cdots+\lambda_k=n$.
A skew Young diagram of shape $\lambda/\mu$ is a collection of row cells such that the $k$-th row begins at the $(\mu_k+1)$-st position and has $\lambda_k-\mu_k$ cells, where $\lambda=(\lambda_1,\ldots,\lambda_k), \mu=(\mu_1,\ldots,\mu_k)$ are partitions and $\lambda_i\ge\mu_i$ for $1\le i\le k$.
A standard skew Young tableau $T$ of shape $\lambda/\mu$ is a filling of the Young diagram of shape $\lambda/\mu$ with integers $1,2,\ldots,|\lambda|-|\mu|$ such that each row and each column is increasing. When $\mu=\emptyset$, a standard skew Young tableau of shape $\lambda$ is called a standard Young tableau. The set of standard Young tableaux of shape $\lambda$ is denoted by $SYT(\lambda)$.

The   jeu de taquin slide  of Sch\"{u}tzenberger \cite{Sz1} can be described as follows, see also  Stanley \cite{Stanley2}.

Given a standard skew Young tableau $T$ of skew shape $\lambda/\mu$, pick an adjacent empty cell $c$ that can be added to the skew diagram $\lambda/\mu$, that is, $c$  shares at least one edge with some cell in $T$, and $\{c\} \cup \lambda\setminus\mu$ is also a skew diagram.
There are two kinds of slide, the inward slide and the outward slide, depending on whether $c$ lies to the upper left or  the lower right of $T$ .

The inward slide begins with that $c$ lies to the upper left of $T$. Slide the number from its neighbouring cells to its right or below into $c$; if $c$ has both neighbours,  then pick the smaller one.   If the cell that just has been emptied has no neighbour to its right or below, then the slide   terminates.  Otherwise, slide a number into that cell according to the same rule as before, and continue in this way until the slide is terminated.
The outward slide begins with that $c$ lies to the lower right of $T$. Slide the number from its neighbouring cells to its left or above into $c$,  picking the larger one if there is a choice.

Denote the inward (resp., outward) jeu de taquin slide by inward (resp., outward) $jdt$ for short.
The sequences of empty cells during the inward (resp., outward) $jdt$ are called the inward (resp., outward) $jdt$ path. The inward (resp., outward) $jdt$ path is a directed path,  beginning at the first empty cell.

It is easy to see that  after either one of the slide transformation, the resulting tableau (with the added empty cell removed) is still a skew  standard Young tableau.

\begin{defi}
Given  $T\in SYT(\lambda)$ with $\lambda\vdash n$, the promotion operation on $T$, denoted by $\partial(T)$, is defined as follows:

Find the cell $(i,j)$ of $T$ that contains $n$, and delete $(i,j)$ to create an empty cell. Apply outward $jdt$ from $(i,j)$  to create an empty cell at the upper left of $\lambda$. (It is necessarily   the cell $(1,1)$.)  Put 0 in the empty cell. Now add 1 to each entry of the current filling of $\lambda$ to obtain $\partial(T)$.
\end{defi}

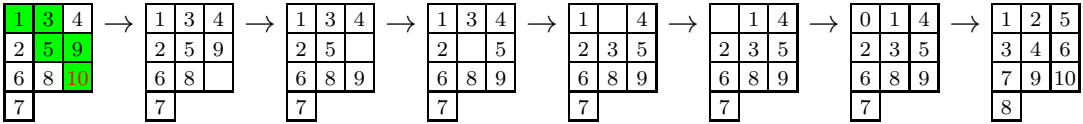
\begin{figure}[!htb]
\setlength{\unitlength}{1.5mm}
\begin{center}
\begin{picture}(100,2)
\ytableausetup{boxsize=.9em}
\ytableausetup{aligntableaux=top}
\begin{ytableau}
*(green)\scriptstyle1 &  *(green)\scriptstyle3 & \scriptstyle4 \\
\scriptstyle2 &  *(green)\scriptstyle5 & *(green)\scriptstyle9 \\
\scriptstyle6 &  \scriptstyle8 & *(green){\color{red}\scriptstyle10}\\
\scriptstyle7
\end{ytableau}
$\rightarrow$
\begin{ytableau}
\scriptstyle1 &  \scriptstyle3 & \scriptstyle4 \\
\scriptstyle2 &  \scriptstyle5 & \scriptstyle9 \\
\scriptstyle6 &  \scriptstyle8 &\empty\\
\scriptstyle7
\end{ytableau}
$\rightarrow$
\begin{ytableau}
\scriptstyle1 &  \scriptstyle3 & \scriptstyle4 \\
\scriptstyle2 &  \scriptstyle5 &   \\
\scriptstyle6 &  \scriptstyle8 & \scriptstyle9 \\
\scriptstyle7
\end{ytableau}
 $\rightarrow$
\begin{ytableau}
\scriptstyle1 &  \scriptstyle3 & \scriptstyle4 \\
\scriptstyle2 &    & \scriptstyle5  \\
\scriptstyle6 &  \scriptstyle8 & \scriptstyle9 \\
\scriptstyle7
\end{ytableau}
$\rightarrow$
\begin{ytableau}
\scriptstyle1 &    & \scriptstyle4 \\
\scriptstyle2 &  \scriptstyle3  & \scriptstyle5  \\
\scriptstyle6 &  \scriptstyle8 & \scriptstyle9 \\
\scriptstyle7
\end{ytableau}
$\rightarrow$
\begin{ytableau}
\empty  &  \scriptstyle1  & \scriptstyle4 \\
\scriptstyle2 &  \scriptstyle3  & \scriptstyle5  \\
\scriptstyle6 &  \scriptstyle8 & \scriptstyle9 \\
\scriptstyle7
\end{ytableau}
$\rightarrow$
\begin{ytableau}
\scriptstyle0  &  \scriptstyle1  & \scriptstyle4 \\
\scriptstyle2 &  \scriptstyle3  & \scriptstyle5  \\
\scriptstyle6 &  \scriptstyle8 & \scriptstyle9 \\
\scriptstyle7
\end{ytableau}
$\rightarrow$
\begin{ytableau}
\scriptstyle1 &  \scriptstyle2  & \scriptstyle5 \\
\scriptstyle3 &  \scriptstyle4  & \scriptstyle6  \\
\scriptstyle7 &  \scriptstyle9 & \scriptstyle10 \\
\scriptstyle8
\end{ytableau}
\end{picture}
\end{center}
\vspace{1cm}
\caption{The construction of the promotion $\partial(T)$.}\label{figp}
\end{figure}

Figure \ref{figp} is an illustration of the construction of the promotion. The outward $jdt$ path  consists of the cells in green.

\begin{defi}
Given  $T\in SYT(\lambda)$ with $\lambda\vdash n$, the dual-promotion operation on $T$, denoted by $\partial^*(T)$, is defined as follows:

Find the cell $(i,j)$ of $T$ that contains $1$, and delete $(i,j)$ to create an empty cell. Apply inward $jdt$ from $(i,j)$ to create an empty cell at the lower right of $\lambda$.  Put $n+1$ in the empty cell. Now subtract 1 to each entry of the current filling of $\lambda$ to obtain $\partial^*(T)$.
\end{defi}

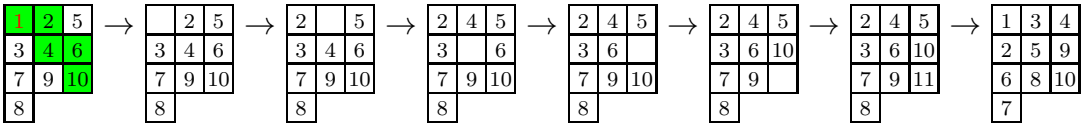
\begin{figure}[!htb]
\setlength{\unitlength}{1.5mm}
\begin{center}
\begin{picture}(100,2)
\ytableausetup{boxsize=.9em}
\ytableausetup{aligntableaux=top}
\begin{ytableau}
*(green){\color{red}\scriptstyle1} &  *(green)\scriptstyle2 & \scriptstyle5 \\
\scriptstyle3 & *(green)\scriptstyle4 & *(green)\scriptstyle6 \\
\scriptstyle7 &  \scriptstyle9 & *(green)\scriptstyle10\\
\scriptstyle8
\end{ytableau}
$\rightarrow$
\begin{ytableau}
\empty &  \scriptstyle2 & \scriptstyle5 \\
\scriptstyle3 & \scriptstyle4 & \scriptstyle6 \\
\scriptstyle7 &  \scriptstyle9 & \scriptstyle10\\
\scriptstyle8
\end{ytableau}
$\rightarrow$
\begin{ytableau}
 \scriptstyle2 &    & \scriptstyle5 \\
\scriptstyle3 & \scriptstyle4 & \scriptstyle6 \\
\scriptstyle7 &  \scriptstyle9 & \scriptstyle10\\
\scriptstyle8
\end{ytableau}
 $\rightarrow$
\begin{ytableau}
 \scriptstyle2 & \scriptstyle4   & \scriptstyle5 \\
\scriptstyle3 &   & \scriptstyle6 \\
\scriptstyle7 &  \scriptstyle9 & \scriptstyle10\\
\scriptstyle8
\end{ytableau}
 $\rightarrow$
\begin{ytableau}
 \scriptstyle2 & \scriptstyle4   & \scriptstyle5 \\
\scriptstyle3 &  \scriptstyle6 &   \\
\scriptstyle7 &  \scriptstyle9 & \scriptstyle10\\
\scriptstyle8
\end{ytableau}
$\rightarrow$
\begin{ytableau}
 \scriptstyle2 & \scriptstyle4   & \scriptstyle5 \\
\scriptstyle3 &  \scriptstyle6 & \scriptstyle10  \\
\scriptstyle7 &  \scriptstyle9 &\\
\scriptstyle8
\end{ytableau}
$\rightarrow$
\begin{ytableau}
 \scriptstyle2 & \scriptstyle4   & \scriptstyle5 \\
\scriptstyle3 &  \scriptstyle6 & \scriptstyle10  \\
\scriptstyle7 &  \scriptstyle9 &\scriptstyle11\\
\scriptstyle8
\end{ytableau}
$\rightarrow$
\begin{ytableau}
 \scriptstyle1 & \scriptstyle3   & \scriptstyle4 \\
\scriptstyle2 &  \scriptstyle5 & \scriptstyle9  \\
\scriptstyle6 &  \scriptstyle8 &\scriptstyle10\\
\scriptstyle7
\end{ytableau}
\end{picture}
\end{center}
\vspace{1cm}
\caption{The construction of the dual-promotion $\partial^*(T)$.}\label{figdp}
\end{figure}

Figure \ref{figdp} is an example of the construction of the dual-promotion. The inward $jdt$ path  consists of the cells in green.

It is easy to see that $\partial$ and $\partial^*$ are inverses to each other.
The following proposition was shown by Edelman and Greene in \cite[Corollary 7.23]{Edelman}.

\begin{prop}\label{pe}
Let $T\in SYT(\delta)$ such that $\delta=(n-1,\ldots,2,1)$. Then
      \[\partial^L(T)=T^t,\]
      where $L=n(n-1)/2$ and $T^t$ is the transpose of $T$. Clearly, $(\partial^*)^L(T)=\partial^L(T)=T^t$.
\end{prop}

We proceed to recall some results in \cite{Edelman}. Given a standard staircase Young tableau $T$ of shape $(n-1,\ldots,2,1)$ and let $L=n(n-1)/2$,  we can use the promotion operator repeatedly to obtain an integer sequence   $a_1,\ldots,a_{L}$ as follows. Let $\partial^{0}(T)=T$. For $1\le k\le L$, at the $k$-th step, apply $\partial$ to $\partial^{k-1}(T)$, denote $a_k$ by the column coordinate of the first cell we deleted, i.e., the cell containing $n$. Let
\begin{align}
\Gamma(T)=s_{a_1}s_{a_2}\cdots s_{a_{L}}.
\end{align}
Edelman and Greene \cite[Theorem 5.4]{Edelman} showed that

\begin{theo} \label{gm}
The map $\Gamma$ is a bijection between $SYT(\delta)$ and $R(w_0)$, where $w_0=n\cdots 21$ and $\delta=(n-1,\ldots,2,1)$.
\end{theo}

In fact, in \cite[Theorem 5.4]{Edelman}, $\Gamma(T)$ is defined to be equal to $s_{a_{L}}\cdots s_{a_2}s_{a_1}$, which is also a reduced word for $w_0$, since $w_0^{-1}=w_0$. However, we find it is more convenient to reverse the order of the sequence $s_{a_{L}}\cdots s_{a_2}s_{a_1}$ for our purpose.

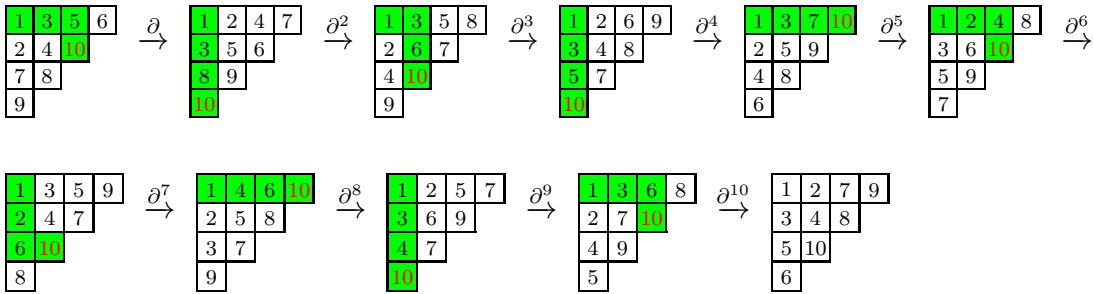
\begin{figure}[!htb]
\setlength{\unitlength}{1.5mm}
\begin{center}
\begin{picture}(140,5)

\put(0,0){
\ytableausetup{boxsize=.85em}
\ytableausetup{aligntableaux=top}
\begin{ytableau}
*(green)\scriptstyle1 &  *(green)\scriptstyle3 & *(green)\scriptstyle5 & \scriptstyle6\\
\scriptstyle2 &  \scriptstyle4 & *(green)\scriptstyle{\color{red}10} \\
\scriptstyle7 &  \scriptstyle8  \\
\scriptstyle9
\end{ytableau}
\
$\underrightarrow{\scriptstyle\partial}$
\
\begin{ytableau}
*(green)\scriptstyle1 &  \scriptstyle2 & \scriptstyle4 & \scriptstyle7\\
*(green)\scriptstyle3 &  \scriptstyle5 & \scriptstyle6 \\
*(green)\scriptstyle8 &  \scriptstyle9  \\
*(green)\scriptstyle{\color{red}10}
\end{ytableau}
\
$\underrightarrow{\scriptstyle\partial^2}$
\
\begin{ytableau}
*(green)\scriptstyle1 &  *(green)\scriptstyle3 & \scriptstyle5 & \scriptstyle8\\
\scriptstyle2 &  *(green)\scriptstyle6 & \scriptstyle7 \\
\scriptstyle4 &  *(green)\scriptstyle{\color{red}10}  \\
\scriptstyle9
\end{ytableau}
\
$\underrightarrow{\scriptstyle\partial^3}$
\
\begin{ytableau}
*(green)\scriptstyle1 &  \scriptstyle2 & \scriptstyle6 & \scriptstyle9\\
*(green)\scriptstyle3 &  \scriptstyle4 & \scriptstyle8 \\
*(green)\scriptstyle5 &  \scriptstyle7  \\
*(green)\scriptstyle{\color{red}10}
\end{ytableau}
\
$\underrightarrow{\scriptstyle\partial^4}$
\
\begin{ytableau}
*(green)\scriptstyle1 &  *(green)\scriptstyle3 & *(green)\scriptstyle7 & *(green)\scriptstyle{\color{red}10}\\
\scriptstyle2 &  \scriptstyle5 & \scriptstyle9 \\
\scriptstyle4 &  \scriptstyle8  \\
\scriptstyle6
\end{ytableau}
\
$\underrightarrow{\scriptstyle\partial^5}$
\
\begin{ytableau}
*(green)\scriptstyle1 &  *(green)\scriptstyle2 & *(green)\scriptstyle4 & \scriptstyle8\\
\scriptstyle3 &  \scriptstyle6 & *(green)\scriptstyle{\color{red}10} \\
\scriptstyle5 &  \scriptstyle9  \\
\scriptstyle7
\end{ytableau}
\
$\underrightarrow{\scriptstyle\partial^6}$
\
}

\put(0,-15){
\ytableausetup{boxsize=.9em}
\ytableausetup{aligntableaux=top}
%
\begin{ytableau}
*(green)\scriptstyle1 &  \scriptstyle3 & \scriptstyle5 & \scriptstyle9\\
*(green)\scriptstyle2 &  \scriptstyle4 & \scriptstyle7 \\
*(green)\scriptstyle6 &  *(green)\scriptstyle{\color{red}10}  \\
\scriptstyle8
\end{ytableau}
\
$\underrightarrow{\scriptstyle\partial^7}$
\
\begin{ytableau}
*(green)\scriptstyle1 &  *(green)\scriptstyle4 & *(green)\scriptstyle6 &*(green)\scriptstyle {\color{red}10}\\
\scriptstyle2 &  \scriptstyle5 & \scriptstyle8 \\
\scriptstyle3 &  \scriptstyle7  \\
\scriptstyle9
\end{ytableau}
\
$\underrightarrow{\scriptstyle\partial^8}$
\
\begin{ytableau}
*(green)\scriptstyle1 &  \scriptstyle2 & \scriptstyle5 & \scriptstyle7\\
*(green)\scriptstyle3 &  \scriptstyle6 & \scriptstyle9 \\
*(green)\scriptstyle4 &  \scriptstyle7  \\
*(green)\scriptstyle{\color{red}10}
\end{ytableau}
\
$\underrightarrow{\scriptstyle\partial^9}$
\
\begin{ytableau}
*(green)\scriptstyle1 &  *(green)\scriptstyle3 & *(green)\scriptstyle6 & \scriptstyle8\\
\scriptstyle2 &  \scriptstyle7 & *(green)\scriptstyle{\color{red}10} \\
\scriptstyle4 &  \scriptstyle9  \\
\scriptstyle5
\end{ytableau}
\
$\underrightarrow{\scriptstyle\partial^{10}}$
\
\begin{ytableau}
\scriptstyle1 &  \scriptstyle2 & \scriptstyle7 & \scriptstyle9\\
\scriptstyle3 &  \scriptstyle4 & \scriptstyle8 \\
\scriptstyle5 &  \scriptstyle10  \\
\scriptstyle6
\end{ytableau}
}
\end{picture}
\end{center}
\vspace{3cm}
\caption{The construction of $\Gamma(T)$.}\label{fig2sf}
\end{figure}

Figure \ref{fig2sf} is an example of the construction of $\Gamma(T)$. The outward $jdt$ paths   are the cells in green.
The cells containing the first deleted entry (namely, 10) at each step of the  promotion are
\begin{align}\label{sq}
(2,3), (4,1), (3,2), (4,1), (1,4), (2,3), (3,2), (1,4), (4,1), (2,3),
\end{align}
 respectively. Thus
 \[\Gamma(T)=s_3s_1s_2s_1s_4s_3s_2s_4s_1s_3,\]
 which is a reduced expression of $w_0=54321$.

\begin{figure}[!htb]
\setlength{\unitlength}{1.5mm}
\begin{center}
\begin{picture}(140,5)

\put(0,0){
\ytableausetup{boxsize=.9em}
\ytableausetup{aligntableaux=top}
\begin{ytableau}
*(green)\scriptstyle1 &  \scriptstyle3 & \scriptstyle5 & \scriptstyle6\\
*(green)\scriptstyle2 &  *(green)\scriptstyle4 & \scriptstyle10 \\
\scriptstyle7 & *(green) \scriptstyle8  \\
\scriptstyle9
\end{ytableau}
\
$\underrightarrow{\scriptstyle\partial^*}$
\
\begin{ytableau}
*(green)\scriptstyle1 &  *(green)\scriptstyle2 & *(green)\scriptstyle4 & *(green)\scriptstyle5\\
\scriptstyle3 &  \scriptstyle7 & \scriptstyle9 \\
\scriptstyle6 &  \scriptstyle{\color{red}10}  \\
\scriptstyle8
\end{ytableau}
\
$\underrightarrow{\scriptstyle(\partial^*)^2}$
\
\begin{ytableau}
*(green)\scriptstyle1 &  \scriptstyle3 & \scriptstyle4 & \scriptstyle{\color{red}10}\\
*(green)\scriptstyle2 &  \scriptstyle6 & \scriptstyle8 \\
*(green)\scriptstyle5 &  \scriptstyle9  \\
*(green)\scriptstyle7
\end{ytableau}
\
$\underrightarrow{\scriptstyle(\partial^*)^3}$
\
\begin{ytableau}
*(green)\scriptstyle1 &  *(green)\scriptstyle2 & *(green)\scriptstyle3 & \scriptstyle9\\
\scriptstyle4 &  \scriptstyle5 & *(green)\scriptstyle7 \\
\scriptstyle6 &  \scriptstyle8  \\
\scriptstyle{\color{red}10}
\end{ytableau}
\
$\underrightarrow{\scriptstyle(\partial^*)^4}$
\
\begin{ytableau}
*(green)\scriptstyle1 &  *(green)\scriptstyle2 & \scriptstyle6 & \scriptstyle8\\
\scriptstyle3 &  *(green)\scriptstyle4 & \scriptstyle{\color{red}10} \\
\scriptstyle5 &  *(green)\scriptstyle7  \\
\scriptstyle9
\end{ytableau}
\
$\underrightarrow{\scriptstyle(\partial^*)^5}$
\
\begin{ytableau}
*(green)\scriptstyle1 &  \scriptstyle3 & \scriptstyle5 & \scriptstyle7\\
*(green)\scriptstyle2 &  \scriptstyle6 & \scriptstyle9 \\
*(green)\scriptstyle4 &  \scriptstyle{\color{red}10}  \\
*(green)\scriptstyle8
\end{ytableau}

}

\put(0,-15){
\
$\underrightarrow{\scriptstyle(\partial^*)^6}$
\
%
\begin{ytableau}
*(green)\scriptstyle1 &  *(green)\scriptstyle2 & *(green)\scriptstyle4 & *(green)\scriptstyle6\\
\scriptstyle3 &  \scriptstyle5 & \scriptstyle8 \\
\scriptstyle7 &  \scriptstyle9 \\
\scriptstyle{\color{red}10}
\end{ytableau}
\
$\underrightarrow{\scriptstyle(\partial^*)^7}$
\
\begin{ytableau}
*(green)\scriptstyle1 &  \scriptstyle3 & \scriptstyle5 &\scriptstyle {\color{red}10}\\
*(green)\scriptstyle2 &  *(green)\scriptstyle4 & *(green)\scriptstyle7 \\
\scriptstyle6 &  \scriptstyle8  \\
\scriptstyle9
\end{ytableau}
\
$\underrightarrow{\scriptstyle(\partial^*)^8}$
\
\begin{ytableau}
*(green)\scriptstyle1 &  *(green)\scriptstyle2 & *(green)\scriptstyle4 & *(green)\scriptstyle9\\
\scriptstyle3 &  \scriptstyle6 & \scriptstyle{\color{red}10} \\
\scriptstyle5 &  \scriptstyle7  \\
\scriptstyle8
\end{ytableau}
\
$\underrightarrow{\scriptstyle(\partial^*)^9}$
\
\begin{ytableau}
*(green)\scriptstyle1 &  \scriptstyle3 & \scriptstyle8 & \scriptstyle{\color{red}10}\\
*(green)\scriptstyle2 &\scriptstyle5 & \scriptstyle9 \\
*(green)\scriptstyle4 &  *(green)\scriptstyle6  \\
\scriptstyle7
\end{ytableau}
\
$\underrightarrow{\scriptstyle(\partial^*)^{10}}$
\
\begin{ytableau}
\scriptstyle1 &  \scriptstyle2 & \scriptstyle7 & \scriptstyle9\\
\scriptstyle3 &  \scriptstyle4 & \scriptstyle8 \\
\scriptstyle5 &  \scriptstyle{\color{red}10}  \\
\scriptstyle6
\end{ytableau}
}
\end{picture}
\end{center}
\vspace{3cm}
\caption{The construction of $\Gamma^*(T)$.}\label{figiv}
\end{figure}
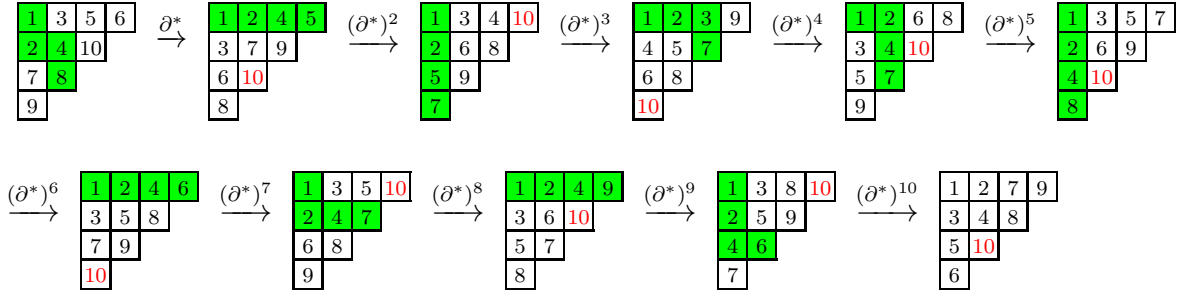

We can obtain the reverse of $\Gamma(T)$ by applying $\partial^*$ to $T$. For example,  Figure \ref{figiv} is an illustration of applying $\partial^*$ ten times to the first tableau in Figure \ref{fig2sf}. The inward $jdt$ paths are the cells in green.

If we record the positions of the last cells of the inward $jdt$ paths at each step,  we would obtain the sequence
\begin{align}\label{seq}
(3, 2), (1, 4), (4, 1), (2, 3), (3, 2), (4, 1), (1, 4), (2, 3), (1, 4), (3, 2).
\end{align}
Compare with \eqref{sq}, we find that the reverse of column coordinates of \eqref{sq}   equals to the row coordinates of the sequences in \eqref{seq}. Let  $\Gamma^*(T)=s_{a_1'}s_{a_2'}\cdots s_{a'_{L}}$, where $a'_1a'_2\cdots a'_{L}$ is the sequence of row coordinates of the last cells at each step of the inward $jdt$ path during applying $\partial^*$ $L$ times to $T$.
 We have

\begin{prop}\label{gama}
Let $T\in SYT(\delta)$ with $\delta=(n-1,\ldots,2,1)$.  Then $\Gamma(T)$ and $\Gamma^*(T)$ are reverses to each other.
\end{prop}

\pf Let $L=n(n-1)/2$.
For $0\le k\le L-1$,
since $\partial^k(T)^t=\partial^L(\partial^k(T))
=\partial^{L+k}(T)$ and
$\partial^{k}(T^t)=\partial^k(\partial^L(T))=\partial^{L+k}(T),$
we have
\[\partial^k(T)^t=\partial^k(T^t).\]
Assume that during applying $\partial$ to $\partial^{k}(T)$ to obtain $\partial^{k+1}(T)$, the first cell we deleted is $(i,j)$.
Then if we apply $\partial^*$ to $\partial^{k+1}(T)$, we would obtain  $\partial^{k}(T)$, and the last cell in the inward $jdt$ path is also $(i,j)$. Clearly, if we we apply $\partial^*$ to $\partial^{k+1}(T)^t$, the last cell in the inward $jdt$ path   is   $(j,i)$. On the other hand,
\[\partial^{k+1}(T)^t=\partial^{k+1}(T^t)
=(\partial^*)^{-k-1}((\partial^*)^L(T))=(\partial^*)^{L-k-1}(T).\]
That is to say, during applying $\partial$ $L$ times to $T$, if the deleted cells are
\[(i_1,j_1),(i_2,j_2),\ldots,(i_L,j_L),\]
then during applying $\partial^*$ $L$ times to $T$, the  last cells in the inward $jdt$ paths are
\[(j_L,i_L),\ldots,(j_{2},i_{2}),(j_1,i_1).\]
This completes the proof.
\qed

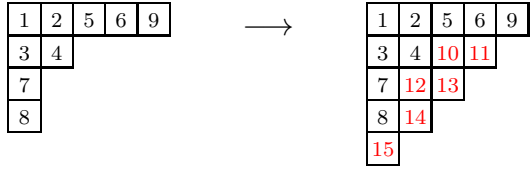
\begin{figure}[!htb]
\setlength{\unitlength}{1.5mm}
\begin{center}
\begin{picture}(50,5)
\ytableausetup{boxsize=1em}
\begin{ytableau}
\scriptstyle1 &  \scriptstyle2 & \scriptstyle5 & \scriptstyle6& \scriptstyle9\\
\scriptstyle3 &  \scriptstyle4   \\
\scriptstyle7   \\
\scriptstyle8
\end{ytableau}
\qquad
$\longrightarrow$
\qquad
\begin{ytableau}
\scriptstyle1 &  \scriptstyle2 & \scriptstyle5 & \scriptstyle6& \scriptstyle9\\
\scriptstyle3 &  \scriptstyle4 & {\color{red}\scriptstyle10}& {\color{red}\scriptstyle11}  \\
\scriptstyle7 & {\color{red}\scriptstyle12}& {\color{red}\scriptstyle13}  \\
\scriptstyle8& {\color{red}\scriptstyle14}\\
{\color{red}\scriptstyle15}
\end{ytableau}
\end{picture}
\end{center}
\vspace{1.3cm}
\caption{A standard Young Tableau $T$ and the corresponding $T^+$.}\label{fig2g}
\end{figure}

For a general Young tableau $T$ with shape $\lambda\vdash n$, in order to obtain a balanced tableau of the same shape as $T$, Edelman and Greene \cite{Edelman} first ``pack" $T$ into a staircase Young tableau $T^+$ in a canonical way. That is, let $\lambda^+$ be the smallest staircase shape which contains $\lambda$. Let $L$ denote the number of cells in $\lambda^+$. Then $T^+$ is the tableau obtained from $T$ by filling the cells of the skew shape $\lambda^+/\lambda$ with $\{n+1,n+2,\ldots,L\}$ from top to bottom and from left to right, see Figure \ref{fig2g} for an example.

The following proposition is a generalization of Theorem \ref{gm}.
 
\begin{prop}\label{gammat}
Let $T\in SYT(\lambda)$ with $\lambda\vdash n$. Assume that $T^+$ has $L$ entries and $\Gamma(T^+)=s_{a_1}s_{a_2}\cdots s_{a_L}$. Then $\Omega(T)=s_{a_{L-n+1}}s_{a_{L-n+2}}\cdots s_{a_L}$ is a reduced word for  the permutation $w=s_{a_{L-n}}\cdots s_{a_2}s_{a_1} w_0$ and $D(w)$ has the same shape as $T$. Moreover, $\Omega$ is a bijection between $SRT(w)$ and $R(w)$.
\end{prop}

\pf Assume that $T^+$ has $N-1$ columns, then $L=N(N-1)/2$ and $w_0=N\cdots 21$. Let $S$ denote the set of cells containing the  added entries $n+1,\ldots,L$ from $T$ to $T^+$. It is easy to see that during applying $\partial$ to $T^+$,  the cells of $S$ always remain in the same column, and the first $|S|$ deleted cells are the cells of $S$ containing the entries $L,\ldots,n+1$, respectively.

More specifically, the indices $a_1,\ldots, a_{L-n}$ of the first $|S|$ entries of $\Gamma(T^+)$ can be obtained by reading the column coordinates of the cells of $S$ in $T^+$ containing the entries $L,\ldots,n+1$, respectively.
For example, if $T$ is   the tableau in Figure \ref{fig2g}, then the first 6 entries of $\Gamma(T^+)$ are $(s_1)(s_2)(s_3s_2)(s_4s_3)$.

Since   $\Gamma(T^+)=s_{a_1}s_{a_2}\cdots s_{a_L}$ is a reduced word for $w_0$, we see that $\Omega(T)=s_{a_{L-n+1}}\cdots$ $ s_{a_{L-1}}s_{a_L}$ is a reduced word for the permutation $w=s_{a_{L-n}}\cdots s_{a_2}s_{a_1} w_0$.
Note that the effect of multiplying $s_{i}$ on the left of a permutation $\pi$ is to interchange the positions of $i$ and $i+1$ in $\pi$.  Thus, if $i+1$ appears to the left of $i$ in $\pi$, the Rothe diagram  $D(s_i\pi)$ can be obtained from $D(\pi)$ by deleting the rightmost cell in the $k$-th row, where $\pi_k=i+1$. In particular, assume that $D(\pi)$ is a Young diagram of shape $\mu=(\mu_1,\ldots,\mu_k)$ such that $\mu_i>\mu_{i+1}$ for some $1\le i\le k$. Then if we delete the last cell in the $i$-th row of $D(\pi)$, we would obtain   $D(s_{\mu_i}\pi)$.

Since $D(w_0)$ is a staircase Young diagram, and $a_1,\ldots, a_{L-n}$  can be obtained by reading the column coordinates of the cells of $S$ in $T^+$ containing $L,\ldots,n+1$, respectively, we see that the Rothe diagram of $w=s_{a_{L-n}}\cdots s_{a_2}s_{a_1} w_0$ can be obtained by deleting the cells of $T^+$ that containing entries $L,\ldots,n+1$ in turn. Consequently,  $\Omega(T)=s_{a_{L-n+1}}\cdots s_{a_{L-1}}s_{a_L}$ is a reduced word for  the permutation $w$ such that $D(w)$ has the same shape as $T$.

It is easy to see that $\Omega$ is an injection from $SRT(w)$ to $R(w)$. On the other hand,
by the results of \cite[Theorem 2.2]{Edelman} and \cite[Theorem 2.4]{Fomin},
we see that
\[|SRT(w)|=|SRT(\lambda)|=|BYT(\lambda)|=|BYT(w)|=|R(w)|.\]
Thus $\Omega$ is a bijection between $R(w)$ and $SRT(w)$, where $D(w)$ is a Young diagram. This completes the proof.
\qed
 
For example, for the standard staircase tableau $T^+$ in Figure \ref{fig2g}, we have $\Gamma(T^+)=s_1s_2s_3s_2s_4s_3s_5s_1s_2s_4s_3s_2s_1s_4s_2$, and then $\Omega(T)=s_5s_1s_2s_4s_3s_2s_1s_4s_2$ is a reduced word for $632415=s_3s_4s_2s_3s_2s_1w_0$.

\section{The equality case}

In this section, we show that if $w$ avoids the patterns 2413, 2431, 3142 and 4132, then $|SRT(w)|$ equals to $|BRT(w)|$, and further give a closed formula for these numbers. Notice that each of these four patterns contains a 132 pattern, thus when $w$ is 132-avoiding, this reduces to result of Edelman and Greene \cite{Edelman} on  the equinumerous between $SYT(\lambda)$ and $BYT(\lambda)$ for any partition $\lambda$.

Recall that the Lehmer code of a permutation $w=w_1\cdots w_n$ is defined by $c(w)=(c_1,\ldots,c_n)$ where
\[c_i=|\{j>i: w_j<w_i\}|.\]
In other words, $c_i$ is the number of cells in the $i$-th row of the Rothe diagram of $w$. Obviously, $c_1+\cdots+c_n=\ell(w)$.
Given two permutations $u=u_1\cdots u_k\in S_k$ and $v=v_1\cdots v_{n-k}\in S_{n-k}$, define the direct sum of $u$ and $v$ to be the permutation
\[u\oplus v=u_1\cdots u_k (v_1+k)\cdots (v_{n-k}+k).\]
 A permutation $w\in S_n$ is called decomposable if it can be expressed as the direct sum of two nontrivial permutations (meaning, $1\le k\le n-1$), and indecomposable otherwise. Every   permutation has a unique expression as the direct sum of indecomposable   permutations. For example, let $w=312486759$.
Then $w$ can be uniquely expressed as $w=w^1\oplus w^2\oplus w^3 \oplus w^4$, where $w^1=312, w^2=1, w^3=4231$ and $w^4=1$.
 
\begin{theo}\label{equals}
Let $w\in S_n$ be a permutation that avoids the patterns 2413, 2431, 3142 and 4132. Assume that $w=w^1\oplus  \cdots \oplus w^k$, where $w^i$ is indecomposable and $c(w^i)=\lambda_i$ for $1\le i\le k$. Then each $\lambda_i$ is a partition and
\[|SRT(w)|=|BRT(w)|={\ell(w) \choose |\lambda_1|,\ldots,|\lambda_k|}\prod_{i=1}^kf^{\lambda_i},\]
where  $f^{\lambda_i}$ is the number of standard Young tableaux of shape $\lambda_i$.
\end{theo}

\pf  We first claim that each $w^i\ (1\le i\le k)$ is 132-avoiding. In fact, since each of $w^1,\ldots,w^k$ is indecomposable,  we need only to show that $w^1$ is 132-avoiding. Suppose to the contrary  that  $w^1$ contains a 132 pattern. We claim that $w^1$ must contain one of the four patterns: 2413, 2431, 3142 or 4132, which is a  contradiction. To prove the claim, let $w^1=w^1_1w^1_2\cdots w^1_m$ and $w^1_{i_1}w^1_{i_2}w^1_{i_3}$ forms a 132 pattern such that $(i_1,i_2,i_3)$ is the lexicographically smallest. There are two cases.

Case 1. $i_1>1$. Since $(i_1,i_2,i_3)$ is the lexicographically smallest, we find that $w^1_{1}>w^1_{i_3}$. Thus there will be a 3142 pattern  (if $w^1_{i_3}<w^1_1<w^1_{i_2}$)  or 4132 pattern  (if $w^1_1>w^1_{i_2}$). The claim follows.

Case 2. $i_1=1$. Since $w^1$ is indecomposable, we have $w^1_1>1$.  Let $w^1_t$ be the largest number appearing before $w^1_{i_2}$. If $w^1_t=w^1_1$, then one of the numbers $1,2,\ldots,w^1_1-1$ must appear after $w^1_{i_2}$, there will be a 2413 or 2431 pattern. Suppose that $w^1_t>w^1_1$. Since $(i_1,i_2,i_3)$ is the lexicographically smallest, we must have $w^1_t<w^1_{i_3}$, and all the numbers $w^1_t+1,\ldots,w^1_{i_3}-1$ appear after $w^1_{i_3}$. Moreover, all the numbers $w^1_1+1,\ldots,w^1_t-1$ appear before $w^1_t$. Since $w^1$ is indecomposable, one of the numbers $1,2,\ldots,w^1_1-1$ must appear after $w^1_{i_2}$. It is easy to check that there will form a 2413 or 2431 pattern again. This proves the claim.

Therefore, each of $c(w^1),\ldots,c(w^k)$ is a partition.
  It is also easy to see that  $D(w)$ can be divided into $k$ nonintersecting Young diagrams (maybe empty), i.e., no two Young diagrams have cells in the same row or column. Since $|SYT(\lambda)|=|BYT(\lambda)|=f^{\lambda}$ for any partition $\lambda$ and $|\lambda_1|+\cdots+|\lambda_k|=\ell(w)$, the theorem follows. \qed
 
It is worth mentioning that in \cite{Atkinson},  Atkinson and   Stitt considered the enumeration of the number of permutations that avoiding 4213, 2413, 3142 and 3241 and obtained a generating function for these numbers. Let $a_n$  denote  the number of permutations in $S_n$ that avoid the four patterns 2413, 2431, 3142 and 4132. After applying a complement and reversing the permutations, we find that $a_n$ equals to the number of permutations avoiding 4213, 2413, 3142 and 3241.
The following proposition is due to \cite{Atkinson}.

\begin{prop} We have
\begin{align*}
\sum_{n=1}^{\infty}a_nx^n&=\frac{-2x}{1-5x+2x^2-(1-x)\sqrt{1-4x}}\\
&=x + 2x^2 + 6x^3 + 20x^4 + 69x^5 + 243x^6 + \cdots.
\end{align*}
\end{prop}

\section{The general case}

In this section, we show that if $w$ contains one of the four patterns 2413, 2431, 3142 and 4132, then $|SRT(w)|$ is strictly smaller than $|BRT(w)|$. The strategy is to establish an injective but not surjective map from $SRT(w)$ to $R(w)$. To this end, we first construct an injection $\eta$ from $SRT(w)$ to $SRT(\widetilde{w})$ such that $\eta(T)\in SRT(\widetilde{w})$ is a Young tableau for any $T\in SRT(w)$, where $\widetilde{w}$ is dominant. Then we can obtain a reduced word of $w$ from $\Omega(\eta(T))\in R(\widetilde{w})$. Finally, we show that this process can not produce all the reduced words of $w$.

In fact, it is enough to show the case that   $w$ is indecomposable. Clearly, if $w$ contains one of the four patterns 2413, 2431, 3142 and 4132, then   $w$   is not a dominant permutation. The idea is to change $w$ into a dominant permutation.
Recall that if $i$ is an ascent of $w$, i.e., $w_i<w_{i+1}$, then $D(ws_i)$ can be obtained from $D(w)$ by first adding a cell $(i,w_{i})$ to $D(w)$, and then moving the cells (if there exist) of $D(w)$ that are in the $(i+1)$-st row and to the right of the column $w_i$ to the $i$-th row.

\noindent{\bf The lifting operation $\eta_i$:}

Let $w=w_1\cdots w_n$ be an indecomposable permutation.
Suppose that $i$ is the first ascent of $w$.
Let $T\in SRT(w)$ be a standard Rothe tableau of $w$.
We  construct a map $\eta_i$ on $T$, called the lifting of $T$ at $(i,w_i)$, as follows.
 \begin{enumerate}
   \item [$(i)$] Apply outward $jdt$  from the empty cell $(i,w_{i})$ of $T$ to create a new empty cell.
       Denote this outward $jdt$ path by $P_{\eta_i}$.
       Since $i$ is the first ascent of $w$, $(c_1,\ldots,c_{i-1})$ forms a partition, the created empty cell is necessarily $(1,1)$.
   \item [$(ii)$] Fill $(1,1)$ with 0, and then add 1 to all the entries of $T$.
   \item [$(iii)$] Move the cells of $T$ that are in the $(i+1)$-st row and to the right of the column $w_i$ to the $i$-th row.
 \end{enumerate}

\begin{figure}[!htb]
\setlength{\unitlength}{1.5mm}
\begin{center}
\begin{picture}(105,5)
\ytableausetup{notabloids}
\ytableausetup{boxsize=.9em}
\begin{ytableau}
*(green)\scriptstyle1 & *(green)\scriptstyle 3 & \scriptstyle6   \\
\scriptstyle2&*(green)\bullet     \\
\scriptstyle4 & \none &\scriptstyle7 & \none &\scriptstyle8 \\
\scriptstyle5
\end{ytableau}
$\rightarrow$
\begin{ytableau}
\scriptstyle1 &  & \scriptstyle 6   \\
\scriptstyle2  & \scriptstyle 3  \\
\scriptstyle4 & \none &\scriptstyle7 & \none &\scriptstyle8 \\
\scriptstyle5
\end{ytableau}
$\rightarrow$
\begin{ytableau}
\empty &  \scriptstyle1  & \scriptstyle6   \\
\scriptstyle2 &\scriptstyle3    \\
\scriptstyle4 & \none &\scriptstyle7 & \none &\scriptstyle8 \\
\scriptstyle5
\end{ytableau}
$\rightarrow$
\begin{ytableau}
\scriptstyle0 &  \scriptstyle1  & \scriptstyle6   \\
\scriptstyle2 &\scriptstyle3    \\
\scriptstyle4 & \none &\scriptstyle7 & \none &\scriptstyle8 \\
\scriptstyle5
\end{ytableau}
$\rightarrow$
\begin{ytableau}
\scriptstyle1 &  \scriptstyle2  & \scriptstyle7   \\
\scriptstyle3 &\scriptstyle4    \\
\scriptstyle5 & \none &\scriptstyle8 & \none &\scriptstyle9 \\
\scriptstyle6
\end{ytableau}
$\rightarrow$
\begin{ytableau}
\scriptstyle1 &  \scriptstyle2  & \scriptstyle7   \\
\scriptstyle3 &\scriptstyle4   &\scriptstyle8 & \none &\scriptstyle9 \\
\scriptstyle5     \\
\scriptstyle6
\end{ytableau}
\end{picture}
\end{center}
\vspace{1cm}
\caption{Applying $\eta_2$ to the tableau Figure \ref{fig1}(b).}\label{fig3ss}
\end{figure}
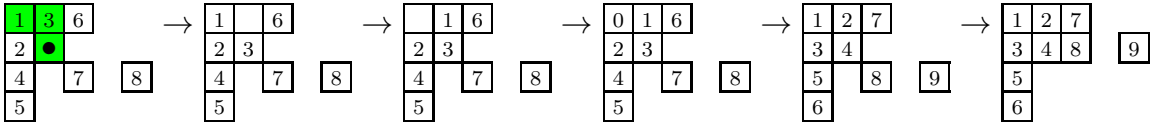

For example,  let $w=426315$. Then   $i=2$ is the first ascent of $w$. Let $T$ be the standard Rothe tableau in Figure \ref{fig1}(b), then $\eta_2(T)$ can be constructed as  in Figure \ref{fig3ss}. We omit the dotted lines of the Rothe diagrams, and the empty cell $(2,w_2)$ of $T$ is   the cell with a dot. The outward $jdt$ path $P_{\eta_2}$ consists of the cells in green.

There are several simple facts about the lifting operation $\eta_i$ that we need.

\noindent{\textbf{Fact 1.}} All the cells $(m,w_i)\notin\eta_i(T)$, where $m\ge i+1$.

\noindent{\textbf{Fact 2.}} In $\eta_i(T)$, the entry in the cell $(i,w_i+1)$ (if exists) is larger than the entries in the cells $(i,w_i)$ and $(i+1,w_i-1)$ (if exists).

Since the lifting operation applies to the first ascent of $w$, it is easy to check that we have the following fact.

\noindent{\textbf{Fact 3.}}
If we can apply $\eta_{i_1}$ to $T$ at $(i_1,j_1)$, and then apply $\eta_{i_2}$ to $\eta_{i_1}(T)$ at $(i_2,j_2)$, and $(i_2,j_2)$ lies to the northeast of $(i_1,j_1)$,  then we must have
$i_1=i_2+1\ \ \text{and}\ \  j_1<j_2.$

\begin{figure}[!htb]
\setlength{\unitlength}{1.5mm}
\begin{center}
\begin{picture}(100,5)
\begin{ytableau}
*(green)\scriptstyle1 &  *(green)\scriptstyle2  & *(green)\scriptstyle7&*(green)\bullet   \\
\scriptstyle3 &\scriptstyle4   &\scriptstyle8 & \none &\scriptstyle9 \\
\scriptstyle5     \\
\scriptstyle6
\end{ytableau}
$\rightarrow$
\begin{ytableau}
\scriptstyle1 &  \scriptstyle2&  & \scriptstyle7   \\
\scriptstyle3 &\scriptstyle4   &\scriptstyle8 & \none &\scriptstyle9 \\
\scriptstyle5     \\
\scriptstyle6
\end{ytableau}
$\rightarrow$
\begin{ytableau}
\scriptstyle1 & & \scriptstyle2   &\scriptstyle7   \\
\scriptstyle3 &  \scriptstyle4   &\scriptstyle8 & \none &\scriptstyle9 \\
\scriptstyle5     \\
\scriptstyle6
\end{ytableau}
$\rightarrow$
\begin{ytableau}
\scriptstyle0 &\scriptstyle1   & \scriptstyle2   &\scriptstyle7   \\
\scriptstyle3 &  \scriptstyle4   &\scriptstyle8 & \none &\scriptstyle9 \\
\scriptstyle5     \\
\scriptstyle6
\end{ytableau}
$\rightarrow$
\begin{ytableau}
\scriptstyle1 &\scriptstyle2   & \scriptstyle3   &\scriptstyle8   \\
\scriptstyle4 &  \scriptstyle5   &\scriptstyle9 & \none &\scriptstyle10 \\
\scriptstyle6     \\
\scriptstyle7
\end{ytableau}
$\rightarrow$
\begin{ytableau}
\scriptstyle1 &\scriptstyle2& \scriptstyle3&\scriptstyle8  &\scriptstyle10   \\
\scriptstyle4 &\scriptstyle5&\scriptstyle9  \\
\scriptstyle6     \\
\scriptstyle7
\end{ytableau}
\end{picture}
\end{center}
\vspace{1cm}
\caption{Applying $\eta_1$ to $\eta_2(T)$.}\label{figs3s}
\end{figure}
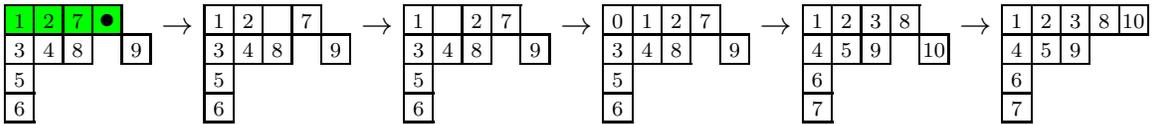

By Fact 2, it is easy to see that $\eta_i(T)\in SRT(ws_i)$, and $\eta_i$ is an injection from $SRT(w)$ to $SRT(ws_i)$.
We can apply the lifting operations $\eta_i$ repeatedly to ``kill'' the ascents of $w$ to change $w$ into a dominant permutation, and change $T$ into a standard Young tableau eventually. The worst situation is that we have to change $w$ into $w_0$, but most of the time we do not need to.
For the running example in Figure \ref{fig3ss}, we need apply $\eta_1$ to $\eta_2(T)$ further to reach a standard Young tableau, see Figure \ref{figs3s}.

For simplicity,   denote $\eta_i(T)=Ts_i$, and let
\begin{align}
\eta(T)=\eta_{i_k}\circ\cdots\circ\eta_{i_1}(T)
=Ts_{i_1}s_{i_2}\cdots s_{i_k},
\end{align}
 where the   $\widetilde{w}=ws_{i_1}s_{i_2}\cdots s_{i_k}$ is a dominant permutation and $Ts_{i_1}s_{i_2}\cdots s_{i_k}\in SRT(\widetilde{w})$ is a standard Young tableau. Obviously, $\eta$ is an injection from $SRT(w)$ to $SRT(\widetilde{w})$.
Now we can state the main theorem of this section.

\begin{theo}\label{inj}
Let $w$ be an indecomposable permutation. Suppose that $T\in SRT(w)$, and   $\eta(T)=Ts_{i_1}s_{i_2}\cdots s_{i_k}$ is a standard Young tableau. Then $\Omega(\eta(T))$ ends with $s_{i_1}s_{i_2}\cdots s_{i_k}$.
\end{theo}

\pf By Proposition \ref{gama} and Proposition \ref{gammat}, it is enough to show that $\Gamma^*(\eta(T))$ begins with $s_{i_k}\cdots s_{i_2} s_{i_1}$. That is, after   ``packing'' $\eta(T)$ into a staircase tableau $\eta(T)^+$,   we apply dual-promotions $\partial^*$ repeatedly to $\eta(T)^+$, we need to show that the first $k$ inward $jdt$ paths end in rows $i_k,\ldots,i_2,i_1$, respectively.

Let $w^0=w$ and $T^0=T$.
For $1\le j\le k$, let  $w^j=ws_{i_1}\cdots s_{i_j}$ and $T^j=Ts_{i_1}\cdots s_{i_j}$.
Recall that $T^j$ is obtained from $T^{j-1}$ by first applying outward $jdt$ beginning at the cell $(i_j,w^{j-1}_{i_j})$, and then moving the cells of $T^{j-1}$ that are in the $(i_j+1)$-st row and to the right of the column $w^{j-1}_{i_j}$ to the $i_j$-th row.

We aim to show that  for $1\le r\le k$, when we apply $\partial^*$ to $(\partial^*)^{r-1}(\eta(T)^+)$,  the inward $jdt$ path $P_r$ must contain the $(k-r+1)$-st outward $jdt$ path $P_{\eta_{i_{k-r+1}}}$, so $P_r$ contains the cell  $(i_{k-r+1},w^{k-r}_{i_{k-r+1}})$. Moreover, after $P_r$ arrives at $(i_{k-r+1},w^{k-r}_{i_{k-r+1}})$, it will not go downward anymore and ends at the $i_{k-r+1}$-th row eventually.

We make induction on $r$. Let us consider the case $r=1$ first.
Denote the set of cells containing the added numbers when we ``pack'' $\eta(T)$ into   $\eta(T)^+$ by $S$. Clearly, all the entries in the cells of $S$ are larger than the entries of $\eta(T)$.
Moreover, by Fact 1,  the cell $(i_k+1,w^{k-1}_{i_k})\in S$, then all the cells $(i_k+1,b)\in S$, where $b\ge w^{k-1}_{i_k}$.
Since the outward $jdt$ and inward $jdt$ are inverses to each other, it is easy to see that the first inward $jdt$ path $P_1$ contains the last outward $jdt$ path $P_{\eta_{i_k}}$, and after $P_1$ arrives at $(i_k,w^{k-1}_{i_k})$, it will end in the $i_k$-th row.

Assume that after applying $\partial^*$  $r\ (r\ge1)$ times to $\eta(T)^+$, the first $r$ inward $jdt$ paths $P_1,\ldots,P_{r}$ contain the outward $jdt$ paths $P_{\eta_{i_k}},\ldots,P_{\eta_{i_{k-r+1}}}$, arrive at the cells
\begin{align}\label{hsa}
(i_k,w^{k-1}_{i_k}),\ldots,(i_{k-r+1},w^{k-r}_{i_{k-r+1}}),
\end{align}
 and end  in rows $i_k,\ldots,$ $i_{k-r+1}$, respectively.
 Firstly, we show that the  $(r+1)$-st inward $jdt$ path $P_{r+1}$
will contain the $(k-r)$-th outward $jdt$ path $P_{\eta_{i_{k-r}}}$ and arrive at $(i_{k-r},w^{k-r-1}_{i_{k-r}})$.

 For $0\le j\le k$, let $\gamma^j$ be the standard Young tableau formed by the connected component of $T^j$ containing $(1,1)$. Clearly, $\gamma^k=\eta(T)$. By the construction of the lifting operation, if we apply $\partial^*$ to $\gamma^j$ ($j\ge 1$), the inward $jdt$ path  will contain the $j$-th outward $jdt$ path $P_{\eta_{i_j}}$, which contains the cell $(i_j,w^{j-1}_{i_j})$. Moreover,  if we delete the cells $(i_j,b)\ (b\ge w^{j-1}_{i_j})$ from $\partial^*(\gamma^j)$, we would obtain $\gamma^{j-1}$.

Therefore,  if we delete  the cells
\begin{align}\label{dt}
\{(i_t,b): k-r+1\le t\le k\ \text{and}\  b\ge w^{t-1}_{i_t}\},
\end{align}
of $(\partial^*)^r(\eta(T)^+)$ and also delete the cells of $S$, we would obtain $\gamma^{k-r}$. Clearly, if we apply $\partial^*$ to $\gamma^{k-r}$, the inward $jdt$ path would coincide with the  $(k-r)$-th outward $jdt$ path $P_{\eta_{i_{k-r}}}$ and arrives at $(i_{k-r},w^{k-r-1}_{i_{k-r}})$. Since $P_{r+1}$ is a directed path and begins at $(1,1)$, it is easy to see that the cells of $S$ have no effect on the route of $P_{r+1}$. Thus   we need to show that the cells  in \eqref{dt}  have no effect on the route of $P_{r+1}$.

In fact, it will be seen that it suffices to show that the cells
\begin{align}\label{xy}
\{(i_t,w^{t-1}_{i_t}): k-r+1\le t\le k \},
\end{align}
have no effect on the route of $P_{r+1}$. There are two situations that one of the cells in \eqref{xy}, say $(i_t,w^{t-1}_{i_t})$, may change the route of $P_{r+1}$.

Case 1. $P_{r+1}$ should go downward, but $(i_t,w^{t-1}_{i_t})$ makes it go right. We claim that this can not happen. It is obvious that $(i_t,w^{t-1}_{i_t})$ lies to north of $(i_{k-r},w^{k-r-1}_{i_{k-r}})$. If $(i_t,w^{t-1}_{i_t})$ lies to the northwest of $(i_{k-r},w^{k-r-1}_{i_{k-r}})$, then $(i_t,w^{t-1}_{i_t})$ does not appear in \eqref{hsa}, and can not affect the route of $P_{r+1}$.
If $(i_t,w^{t-1}_{i_t})$ lies to the northeast of $(i_{k-r},w^{k-r-1}_{i_{k-r}})$, then we must have
$w^{t-1}_{i_t}=w^{k-r-1}_{i_{k-r}}+1$.
By Fact 3, we have
\[i_t=i_{k-r}-1\  \ \text{and}\ \ w^{t-1}_{i_t}=w^{k-r-1}_{i_{k-r}}+1.\]
Since the lifting operation applies to the first ascent of $w$, we see that $(i_t,w^{t-1}_{i_t})$ is one of the cells in \eqref{hsa}. Thus the   $(k-t+1)$-st inward $jdt$ path $P_{k-t+1}$ appears before $P_{r+1}$.
By Fact 2, the cell $(i_t,w^{t-1}_{i_t}+1)$ of $(\partial^*)^{k-t}(\eta(T)^+)$ was either moved above from the $(i_t+1)$-st row of $T^{t-1}$ or belongs to $S$, we see that
the entry in  $(i_t,w^{t-1}_{i_t}+1)$  is larger than the entry in $(i_t+1,w^{t-1}_{i_t}-1)$ in $(\partial^*)^{k-t}(\eta(T)^+)$. By the induction hypothesis, $P_{k-t+1}$ contains  $(i_t,w^{t-1}_{i_t})$ and ends in the $i_t$-th row. Thus, in $(\partial^*)^{k-t+1}(\eta(T)^+)$, the entry in  $(i_t,w^{t-1}_{i_t})$ (which equals to the entry in $(i_t,w^{t-1}_{i_t}+1)$ of $(\partial^*)^{k-t}(\eta(T)^+)$ minus one),  is larger than the entry in $(i_t+1,w^{t-1}_{i_t}-1)$ (which equals to the entry in $(i_t+1,w^{t-1}_{i_t}-1)$ of $(\partial^*)^{k-t}(\eta(T)^+)$ minus one).
Therefore, after $P_{r+1}$ arrives at $(i_t,w^{t-1}_{i_t}-1)$, it will not go right. This proves the claim.

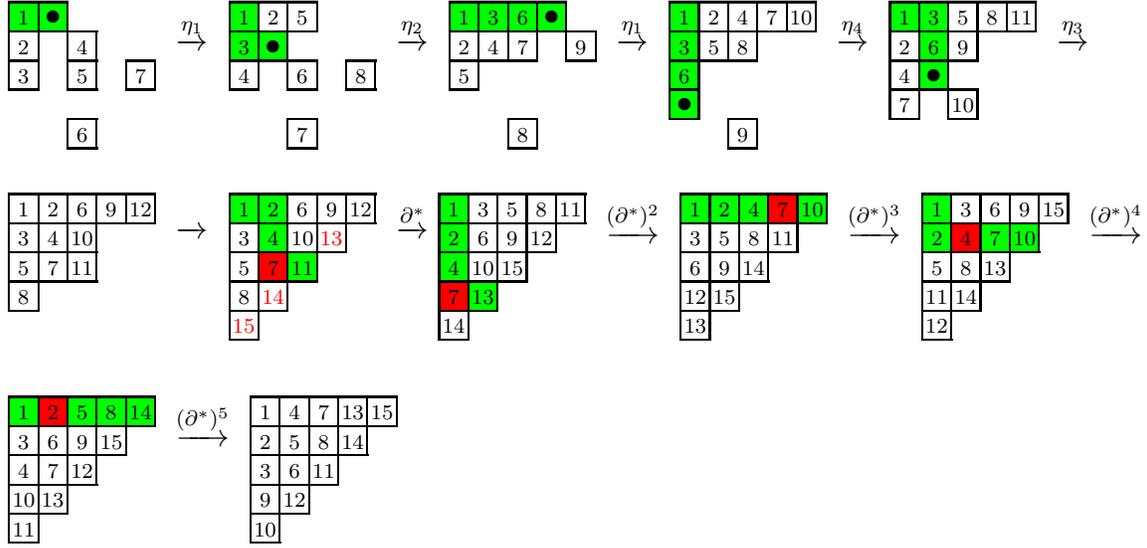
\begin{figure}[!htb]
\setlength{\unitlength}{1.5mm}
\begin{center}
\begin{picture}(90,35)
\ytableausetup{boxsize=.9em}
\ytableausetup{aligntableaux=top}
\put(-5,30){
\begin{ytableau}
*(green)\scriptstyle1 &*(green)\bullet    \\
\scriptstyle2 &\none &\scriptstyle4  \\
\scriptstyle3 &\none &\scriptstyle5&\none &\scriptstyle7\\
\none\\
\none&\none  &\scriptstyle6
\end{ytableau}
$\ \underrightarrow{\scriptstyle\eta_1}\ $
\begin{ytableau}
*(green)\scriptstyle1 & \scriptstyle2&\scriptstyle5    \\
*(green)\scriptstyle3 &*(green)\bullet   \\
\scriptstyle4 &\none &\scriptstyle6&\none &\scriptstyle8\\
\none\\
\none&\none  &\scriptstyle7
\end{ytableau}
$\ \underrightarrow{\scriptstyle\eta_2}\ $
\begin{ytableau}
*(green)\scriptstyle1 &*(green)\scriptstyle3&*(green)\scriptstyle6 &*(green) \bullet   \\
\scriptstyle2 & \scriptstyle4&\scriptstyle7&\none &\scriptstyle9   \\
\scriptstyle5\\
\none\\
\none&\none &\scriptstyle8
\end{ytableau}
$\ \underrightarrow{\scriptstyle\eta_1}\ $
\begin{ytableau}
*(green)\scriptstyle1 &\scriptstyle2&\scriptstyle4 & \scriptstyle7& \scriptstyle10  \\
*(green)\scriptstyle3 & \scriptstyle5&\scriptstyle8   \\
*(green)\scriptstyle6\\
 *(green)\bullet\\
\none&\none  &\scriptstyle9
\end{ytableau}
$\ \underrightarrow{\scriptstyle\eta_4}\ $
\begin{ytableau}
*(green)\scriptstyle1 &*(green)\scriptstyle3&\scriptstyle5 & \scriptstyle8& \scriptstyle11  \\
\scriptstyle2 & *(green)\scriptstyle6&\scriptstyle9   \\
\scriptstyle4&*(green)\bullet\\
\scriptstyle7&\none  &\scriptstyle10
\end{ytableau}
$\ \underrightarrow{\scriptstyle\eta_3}\ $

}
\put(-5,13){
\begin{ytableau}
\scriptstyle1 &\scriptstyle2&\scriptstyle6 & \scriptstyle9& \scriptstyle12  \\
\scriptstyle3 & \scriptstyle4&\scriptstyle10   \\
\scriptstyle5&\scriptstyle7&\scriptstyle11\\
\scriptstyle8
\end{ytableau}
$\ \underrightarrow{ }\ $
\begin{ytableau}
*(green)\scriptstyle1 &*(green)\scriptstyle2&\scriptstyle6 & \scriptstyle9& \scriptstyle12  \\
\scriptstyle3 & *(green)\scriptstyle4&\scriptstyle10 &{\color{red}\scriptstyle13}  \\
\scriptstyle5&*(red)\scriptstyle7&*(green)\scriptstyle11\\
\scriptstyle8&{\color{red}\scriptstyle14}\\
{\color{red}\scriptstyle15}
\end{ytableau}
$\ \underrightarrow{\scriptstyle\partial^*}$
\begin{ytableau}
*(green)\scriptstyle1 &\scriptstyle3 &\scriptstyle5 &\scriptstyle8 &\scriptstyle11 \\
*(green)\scriptstyle2 &\scriptstyle6 &\scriptstyle9 &{\scriptstyle12} \\
*(green)\scriptstyle4&{\scriptstyle10}&{\scriptstyle15}\\
*(red){\scriptstyle7}&*(green){\scriptstyle13}\\
\scriptstyle14
\end{ytableau}
$\ \underrightarrow{\scriptstyle(\partial^*)^2}\ $
\begin{ytableau}
*(green)\scriptstyle1 &*(green)\scriptstyle2 &*(green)\scriptstyle4 &*(red)\scriptstyle7 &*(green)\scriptstyle10 \\
\scriptstyle3 &\scriptstyle5 &\scriptstyle8 &{\scriptstyle11} \\
 \scriptstyle6&{\scriptstyle9}&{\scriptstyle14}\\
{\scriptstyle12}&{\scriptstyle15}\\
{\scriptstyle13}
\end{ytableau}
$\ \underrightarrow{\scriptstyle(\partial^*)^3}\ $
\begin{ytableau}
*(green)\scriptstyle1 &\scriptstyle3 &\scriptstyle6 &\scriptstyle9 &\scriptstyle15 \\
*(green)\scriptstyle2 &*(red)\scriptstyle4 &*(green)\scriptstyle7 &*(green)\scriptstyle10 \\
\scriptstyle5&{\scriptstyle8}&{\scriptstyle13}\\
{\scriptstyle11}&{\scriptstyle14}\\
{\scriptstyle12}
\end{ytableau}
$\ \underrightarrow{\scriptstyle(\partial^*)^4}\ $
}

\put(-5,-5){
\begin{ytableau}
*(green)\scriptstyle1 &*(red)\scriptstyle2 &*(green)\scriptstyle5 &*(green)\scriptstyle8 &*(green)\scriptstyle14 \\
\scriptstyle3 &\scriptstyle6 &\scriptstyle9 &{\scriptstyle15} \\
\scriptstyle4&{\scriptstyle7}&{\scriptstyle12}\\
{\scriptstyle10}&{\scriptstyle13}\\
{\scriptstyle11}
\end{ytableau}
$\ \underrightarrow{\scriptstyle(\partial^*)^5}\ $
\begin{ytableau}
\scriptstyle1 &\scriptstyle4 &\scriptstyle7 &\scriptstyle13 &\scriptstyle15 \\
\scriptstyle2 &\scriptstyle5 &\scriptstyle8 &{\scriptstyle14} \\
\scriptstyle3&{\scriptstyle6}&{\scriptstyle11}\\
{\scriptstyle9}&{\scriptstyle12}\\
{\scriptstyle10}
\end{ytableau}
}

\end{picture}
\end{center}
\vspace{2.2cm}
\caption{ An illustration of the proof of Theorem \ref{inj}.}\label{figlast}
\end{figure}

For example, let $r=1$ and $k=5$. In Figure \ref{figlast}, $\partial^*(\eta(T)^+)$ is the third tableau in the second row. The cell $(3,2)$ belongs to $P_1$, if we delete the cells $(3,2),(3,3)$ and the cells of $S$ in $\partial^*(\eta(T)^+)$, we would obtain $\gamma^4$, which is the fifth tableau in the first row.
If we apply $\partial^*$ to $\gamma^4$, the inward $jdt$ path will coincide with the fourth outward $jdt$ path, which consists of the green cells of the fourth tableau in the first row. However, the cell $(3,2)$ would not affect the route of $P_2$, since in $\partial^*(\eta(T)^+)$, the entry in the cell $(3,2)$ was moved above from the fourth row of $T^4$ (the fifth tableau in the first row), which is larger than the entry in $(4,1)$ of $\partial^*(\eta(T)^+)$.

Case 2. $P_{r+1}$ should go right, but $(i_t,w^{t-1}_{i_t})$ makes it go downward. Clearly, $(i_t,w^{t-1}_{i_t})$ lies to the west of $(i_{k-r},w^{k-r-1}_{i_{k-r}})$. If $(i_t,w^{t-1}_{i_t})$ lies to the northwest of $(i_{k-r},w^{k-r-1}_{i_{k-r}})$, then $(i_t,w^{t-1}_{i_t})$ does not appear in \eqref{hsa}, and can not affect the route of $P_{r+1}$. If $(i_t,w^{t-1}_{i_t})$ lies to the southwest of $(i_{k-r},w^{k-r-1}_{i_{k-r}})$, then by Fact 3  we have
\[i_t=i_{k-r}+1\ \ \text{and}\ \  w^{t-1}_{i_t}<w^{k-r-1}_{i_{k-r}}.\]
By the construction of the lifting operation, it is easy to see that  $(i_t,w^{t-1}_{i_t})$ does not appear in \eqref{hsa} either.

For example, in the second row of Figure \ref{figlast}, the green cells of the fourth tableau consist of the third inward $jdt$ $P_3$. The route of $P_3$ will not affected by the cell $(2,2)$, since $(2,2)$ appears after $(4,1)$ in the inward $jdt$ process.

Next we show that after $P_{r+1}$ arrives at $(i_{k-r},w^{k-r-1}_{i_{k-r}})$, it will end in the $i_{k-r}$-th row at last. There are two cases to consider according to whether the cell $(i_{k-r}+1,w^{k-r-1}_{i_{k-r}})$ belongs to $S$ or not in $\eta(T)^+$.

Case 1. $(i_{k-r}+1,w^{k-r-1}_{i_{k-r}})\in S$. It is obvious that all the cells to the right of $(i_{k-r}+1,w^{k-r-1}_{i_{k-r}})$ also belong to $S$. It is easy to see that $P_{r+1}$ will end in the $i_{k-r}$-th row in this case.

Case 2. $(i_{k-r}+1,w^{k-r-1}_{i_{k-r}})\notin S$. We claim that $(i_{k-r}+1,w^{k-r-1}_{i_{k-r}})$ appears in \eqref{hsa}. By Fact 1,  after applying outward $jdt$ beginning at the position $(i_{k-r},w^{k-r-1}_{i_{k-r}})$, the cell $(i_{k-r}+1,w^{k-r-1}_{i_{k-r}})$ is empty. Moreover, the only way that $(i_{k-r}+1,w^{k-r-1}_{i_{k-r}})$ can belong to $\eta(T)$ is to apply an outward $jdt$ beginning at this cell. Therefore,  if   $(i_{k-r}+1,w^{k-r-1}_{i_{k-r}})\notin S$ in $\eta(T)^+$, this means that we have applied an outward $jdt$ beginning at this cell.
Thus,  $(i_{k-r}+1,w^{k-r-1}_{i_{k-r}})$ appears in \eqref{hsa}.
Let $i_{k-r}+1 =i_{k-m+1}$ for some integer $m$ such that $m\le r$. Then the $m$-th inward $jdt$ path  $P_m$  contains $(i_{k-r}+1,w^{k-r-1}_{i_{k-r}})$ and appears before $P_{r+1}$.

Let $(i_{k-r},b)$  be a cell to the right of $(i_{k-r},w^{k-r-1}_{i_{k-r}})$ in $(\partial^*)^{r}(\eta(T)^+)$, where $b\ge w^{k-r-1}_{i_{k-r}}$. It suffices to show that, in  $(\partial^*)^{r}(\eta(T)^+)$, the entry  in the cell $(i_{k-r},b+1)$ is smaller than the entry in the cell $(i_{k-r}+1,b)$, if exists. Since    $P_m$ appears before $P_{r+1}$, and by the induction hypothesis, $P_m$ ends at the $(i_{k-r}+1)$-st row, we see that the cell   $(i_{k-r}+1,b)$ of $(\partial^*)^{r}(\eta(T)^+)$ belongs to $P_m$ and appears at the position $(i_{k-r}+1,b+1)$ of $(\partial^*)^{m-1}(\eta(T)^+)$.
 Of course, in $(\partial^*)^{m-1}(\eta(T)^+)$, the entry in the cell $(i_{k-r},b+1)$ is smaller than the entry in the cell $(i_{k-r}+1,b+1)$.    This completes the proof. \qed

 Figure \ref{figlast} gives an example of the proof of Theorem \ref{inj}. The first  tableau $T$ belongs to $SRT(w)$ for $w=246153$. After applying $\eta_1,\eta_2,\eta_1,\eta_4,\eta_3$ in turn to $T$, we obtain the standard Young tableau $\eta(T)$ as the first tableau in the second row of Figure \ref{figlast}. The  cells $(1,2),(2,2),(1,4),(4,1),(3,2)$ containing a dot are the cells of the outward $jdt$ beginning with at each step of the lifting operation. When we apply $\partial^*$ to $\eta(T)^+$ repeatedly, the 5 inward $jdt$ paths  arrive at the cells $(3,2),(4,1),(1,4),(2,2),(1,2)$ and end in the rows $3,4,1,2,1$, respectively.

By Theorem \ref{inj}, it is obvious that
$\eta$ is an injection from $SRT(w)$ to
\begin{align}\label{omega}
\{M\in SRT(\widetilde{w}):\Omega(M)\ \text{ends with}\ s_{i_1}s_{i_2}\cdots s_{i_k}\}.
\end{align}
 Since there is a bijection $\zeta$ between $R(w)$ and
\begin{align}\label{pi}
\{\tau\in R(\widetilde{w}): \tau \ \text{ends with}\ s_{i_1}s_{i_2}\cdots s_{i_k}\},
\end{align}
and $\Omega$ is a bijection between the two sets \eqref{omega} and \eqref{pi},
we see that
\begin{align}
|R(w)|=|\{M\in SRT(\widetilde{w}):\Omega(M)\ \text{ends with}\ s_{i_1}s_{i_2}\cdots s_{i_k}\}|
\end{align}
and
\begin{align}
\zeta\circ\Omega\circ\eta: SRT(w)\xlongrightarrow{\eta} SRT(\widetilde{w})\xlongrightarrow{\Omega} R(\widetilde{w})\xlongrightarrow{\zeta} R(w)
\end{align}
is an injection from $SRT(w)$ to $R(w)$.  Thus
\begin{align}
|SRT(w)|\le|R(w)|=|BRT(w)|.
\end{align}
The following theorem, combining with Theorem \ref{equals}, shows that $|SRT(w)|=|BRT(w)|$ if and only if $w$ avoids the four patterns 2431, 2413, 3142  and 4132.
 
\begin{theo}
Let $w$ be a permutation that contains one of the four patterns 2413, 2431, 3142 and 4132. Then $|SRT(w)|$ is strictly smaller than $|BRT(w)|$.
\end{theo}

\pf It is enough to consider the case that $w$ is indecomposable. Given $T\in SRT(w)$, suppose that $\eta(T)=Ts_{i_1}\cdots s_{i_k}$ is a standard Young tableau of $SRT(ws_{i_1}\cdots s_{i_k})$, and $\Omega(\eta(T))$ ends with $s_{i_1}\cdots s_{i_k}$. In fact, it suffices to show the case that $k=1$, i.e., $\eta(T)=Ts_{i}$ is a standard Young tableau of $SRT(ws_{i})$ and $\Omega(\eta(T))$ ends with $s_{i}$.
That is, we need to show that
\begin{align*}
|SRT(w)|&=|\{M\in SRT(ws_i): M=\eta(T)\ \text{ for some} \ T\in SRT(w)\}|\\
&<|\{M\in SRT(ws_i):\Omega(M)\ \text{ends with}\ s_{i}\}|=|R(w)|=|BRT(w)|,
\end{align*}
where $ws_i$ is a  dominant permutation.

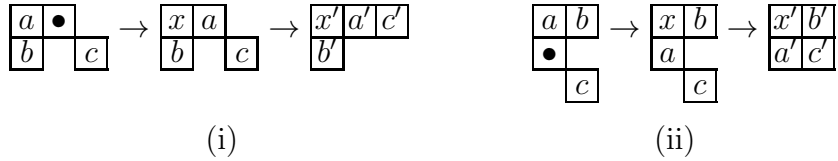
\begin{figure}[!htb]
\setlength{\unitlength}{1.5mm}
\begin{center}
\begin{picture}(70,5)
\ytableausetup{boxsize=1em}
\begin{ytableau}
a & \bullet\\
b&\none &  c   \\
\end{ytableau}
$\rightarrow$
\begin{ytableau}
x&a \\
b&\none &  c
\end{ytableau}
$\rightarrow$
\begin{ytableau}
x' &a' &  c' \\
b'
\end{ytableau}

\put(-18,-10){(i)}

\qquad\qquad

\begin{ytableau}
a &b\\
\bullet  \\
\none &  c   \\
\end{ytableau}
$\rightarrow$
\begin{ytableau}
x &b\\
a\\
\none &  c
\end{ytableau}
$\rightarrow$
\begin{ytableau}
x' &b'\\
a'&c'
\end{ytableau}
\put(-16,-10){(ii)}
\end{picture}
\end{center}
\vspace{1cm}
\caption{The two cases of the lifting operation $\eta_i$.}\label{figeta}
\end{figure}

Recall that in the construction of $\eta(T)=\eta_i(T)$, we first apply an outward $jdt$   from the  empty cell $(i,w_i)$ of $T$ to create the empty cell $(1,1)$. Then we put 0 to $(1,1)$ and add 1 to all the entries of $T$. Finally, we move the cells of $T$ that in the $(i+1)$-st row and to the right of the column $w_i$   to the $i$-th row. There are two cases of the position of the cell $(i,w_i)$  as illustrated in Figure \ref{figeta}. It is possible that a tableau satisfies both cases, see, e.g., the last tableau in the first row of Figure \ref{figlast}.

We refer $a,b,c,x$ to the cells as well as  the entries in the corresponding cells. In both cases in Figure \ref{figeta}, we have $a'=a+1, b'=b+1$, $c'=c+1$  and $x'=x+1$.
The key observation is that since $Ts_i$ is a Young tableau, all the cells $a,b,c$ belong to $T$. Since $T$ is standard, in both cases, we have $c\ge a+2$, then $c'\ge a'+2$. By the proof of Theorem \ref{inj}, when we apply the dual-promotion $\partial^*$ to $\eta(T)^+$, the inward $jdt$ path will contain $x',a',c'$ and ends in the same row as $a'$.

Therefore, the common feature of the tableaux $\eta(T)$ for any $T\in SRT(w)$ is that when we apply $\partial^*$ to $\eta(T)^+$ the inward $jdt$ path contains two adjacent cells $a'$ and $c'$ in the same row and $c'\ge a'+2$.
It is easy to construct a standard Young tableau $M$ of shape $D(ws_i)$ such that when we apply $\partial^*$ to $M^+$, the inward $jdt$ path  contains two adjacent cells $a'$ and $c'$ in the same row but $c'=a'+1$.
This completes the proof. \qed

\noindent{\bf Acknowledgements.}
This work was supported by  the Research Fund for
the Doctoral Program of Higher Education (Grant No. 20130181120103)  and the National Natural Science Foundation of China (Grant No. 11401406).

\end{document}